\documentclass[11pt,twoside]{amsart}

\usepackage{latexsym}
\usepackage{amssymb}

\usepackage{amsfonts,amsmath,amsthm}

\newcommand{\Id}{\mathrm{Id}} 
\newcommand{\eins}{\mathbf{1}}
\newcommand{\dopu}{{:}\allowbreak\ }
\newcommand{\eps}{\varepsilon}
\newcommand{\de}{\delta}
\newcommand{\cal}{\mathcal}

\newcommand{\Om}{(\Omega, \Sigma, \mu)}
\def\DP{Daugavet property}

\newcommand{\loglike}[1]{\mathop{\rm #1}\nolimits}

\newcommand{\dist}{\loglike{dist}}
\newcommand{\lin}{\loglike{lin}}
\newcommand{\supp}{\loglike{supp}}
\newcommand{\Lin}{\loglike{\overline{lin}}}
\newcommand{\conv}{\loglike{conv}}

\newcommand{\re}{\loglike{Re}}

\newcommand{\Daug}{\loglike{Daug}}

\newcommand{\N}{{\mathbb N}}

\newcommand{\calU}{\mathcal{U}}

\theoremstyle{plain}
\newtheorem{thm}{Theorem}[section]

\newtheorem{cor}[thm]{Corollary}
\newtheorem{lemma}[thm]{Lemma}

\theoremstyle{definition}
\newtheorem{definition}[thm]{Definition}

\theoremstyle{remark}

\newtheorem{rem}[thm]{Remark}

\numberwithin{equation}{section}

\newcommand{\begsta}{\begin{statements}}
\newcommand{\begaeq}{\begin{aequivalenz}}
\def\endsta{\end{statements}}
\def\endaeq{\end{aequivalenz}}
\newcommand{\bea}{\begin{eqnarray*}}
\newcommand{\eea}{\end{eqnarray*}}
\newcommand{\kref}[1]{(\ref{#1})}

\newcounter{abc}   
\newcounter{iiiii} 

\newenvironment{aequivalenz}
{\setcounter{iiiii}{0}
\begin{list}%
{{\rm (\roman{iiiii})}}
{\usecounter{iiiii}
\parsep=0pt plus 1pt
\topsep=1pt plus 2pt minus 1pt
\itemsep=1pt plus 2pt minus 1pt
\leftmargin=3\baselineskip
\labelsep=.6\baselineskip
\labelwidth=2.4\baselineskip
\rightmargin 0pt}%
}%
{\end{list}}

\newenvironment{statements}%
{\setcounter{abc}{0}
\begin{list}%
{{\rm (\alph{abc})}}
{\usecounter{abc}
\parsep=0pt plus 1pt
\topsep=1pt plus 2pt minus 1pt
\itemsep=1pt plus 2pt minus 1pt
\leftmargin=3\baselineskip
\labelsep=.6\baselineskip
\labelwidth=2.4\baselineskip
\rightmargin 0pt}%
}%
{\end{list}}

\begin{document}

\title {A Banach space with the Schur and the Daugavet property}

\author{Vladimir Kadets and Dirk Werner}


\address{Faculty of Mechanics and Mathematics, Kharkov National
University,\linebreak
 pl.~Svobody~4,  61077~Kharkov, Ukraine}
\email{vova1kadets@yahoo.com}

\curraddr{Department of Mathematics, Freie Universit\"at Berlin,
Arnim\-allee~\mbox{2--6}, D-14\,195~Berlin, Germany}
\email{kadets@math.fu-berlin.de}

\address{Department of Mathematics, Freie Universit\"at Berlin,
Arnimallee~2--6, \qquad {}\linebreak D-14\,195~Berlin, Germany}
\email{werner@math.fu-berlin.de}

\thanks{The work of the first-named author
was supported by a fellowship from the {\it Alexander-von-Humboldt
Stiftung}.}



\begin{abstract}
We show that a minor refinement of the Bourgain-Rosen\-thal
construction of a Banach space without the Radon-Nikod\'ym property
which contains no bounded $\delta$-trees yields a space with the Daugavet
property and the Schur property. 
Using this example we answer  some open questions
on the structure of such spaces;
in particular we show that the Daugavet
property is not inherited by ultraproducts.
\end{abstract}

\maketitle

\thispagestyle{empty}


\section{Introduction}

A Banach space $X$ is said to have the \textit{Daugavet
property} if every
rank-$1$ operator $T\dopu X \to X$ satisfies
\begin{equation}\label{eq1.1}
\|\Id+T\| =1+ \|T\|.
\end{equation}
Examples include $C[0,1]$, $L_{1}[0,1]$, certain function algebras
such as the disk algebra or $H^\infty$ and some of their subspaces
and quotients. Such spaces are studied in
detail in \cite{KadSSW}.

It is known that a space with the Daugavet
property cannot have the Radon-Nikod\'ym property \cite{Woj92},
that a  space with the Daugavet
property contains a copy of $\ell_1$ \cite{KadSSW} and 
that in a  space with the Daugavet
property not only rank-1 operators, but also weakly compact
operators, strong Radon-Nikod\'ym operators, operators not fixing a copy of
$\ell_1$ and their linear combinations satisfy the Daugavet equation
(\ref{eq1.1}); see \cite{KadSSW}, \cite{KSW}, \cite{Shv1}. A
space with the Daugavet property does not have an unconditional basis
\cite{Kadets} and does not even embed 
into a space with an unconditional basis \cite{KadSSW}.

A Banach space has the Schur property if every weakly convergent
sequence converges in norm. It has been asked in \cite{KKW} and
\cite{dirk-irbull}
whether there exists a Banach space that has both the Schur and the
\DP. Given the isomorphic properties the \DP\ entails as listed
above, the two properties might appear to be mutually exclusive since
a space with the \DP\ should be thought of as rather large
whereas a Schur space, which is hereditarily $\ell_{1}$, could be
imagined as rather thin. In the other direction, previous studies of
the \DP\ \cite{KSW} have supported the conjecture that a space with
the \DP\ always contains a copy of $\ell_{2}$. (This conjecture is
clearly not compatible with the Schur space conjecture.)

In this paper we point out that a certain subspace of $L_{1}$ that
was constructed by Bourgain and Rosenthal \cite{BourRos} indeed has
the Schur and the \DP. This will be shown in Section~2 where we
review the Bourgain-Rosenthal construction along with some minor
modifications to achieve our result. 

In \cite{BKSW} we studied the \DP\ for ultraproducts and asked
whether it is stable under ultraproducts. In Section~3 we show that
direct sums of certain Bourgain-Rosenthal spaces serve as 
counterexamples. To explain why the Bourgain-Rosenthal spaces might
be relevant to this topic we need some technical notions; so we
resume this discussion only in Section~3. 

Throughout the paper we will need the following geometrical
characterisation of the \DP\ from \cite[Lemma~2.2]{KadSSW}.
We denote by $B(X)$ the closed unit ball of $X$ and by $S(X)$ its
unit sphere, and we
recall that a slice of the unit ball is a nonvoid set of the form
$S= \{x\in B(X)\dopu \re x^*(x) > \alpha\}$
for some functional $x^*\in X^*$.

\begin{lemma} \label{lem1.1}
The following assertions are equivalent:
\begaeq
\item
$X$ has the Daugavet property.
\item
For every $x\in S(X)$, $\eps>0$ and every slice $S$
of $B(X)$ there exists some $y\in S$ such that $\|x+y\| \ge 2-\eps$.
\endaeq
\end{lemma}

For an element $x\in S(X)$ and $\eps>0$ denote 
$$ 
l^+(x,\eps) = \{ y\in S(X)\dopu \|x+y\| \ge 2-\eps \}.
$$ 
It is immediate
from Lemma~\ref{lem1.1} and the Hahn-Banach theorem that the \DP\
can be characterised by means of the sets $l^+(x,\eps)$ as follows.
  
\begin{lemma} \label{lem1.2}
The following assertions are equivalent:
\begaeq
\item
$X$ has the Daugavet property.
\item
For every $x\in S(X)$ and $\eps>0$ the set $\conv(l^+(x,\eps))$, the
convex hull of $l^+(x,\eps)$, is dense in $B(X)$.
\endaeq
\end{lemma}


\section{The Bourgain-Rosenthal example revisited}

In this section we will present an example of a Banach space that has
both the Schur and the Daugavet property. It belongs to the class of 
spaces that Bourgain and Rosenthal have constructed in \cite{BourRos}.
However, we base our exposition on the one Benyamini and
Lindenstrauss have given in \cite{BL1} and not on the original paper.

Below we consider the space $L_{1}=L_{1}(\Omega, \Sigma, \mu)$ over a
separable nonatomic measure space. The symbol $\|\,\,.\,\,\|$ will
refer to the $L_{1}$-norm. Besides the norm topology we will also
consider the topology of convergence in measure, generated by the
metric
$$    
d(f,g) = \inf \bigl\{ \eps > 0\dopu  \mu \{t\dopu  |f(t) - g(t)| \ge \eps \} 
  \le   \eps \bigr\}.
$$

We first observe a simple lemma on $L_{1}$-orthogonality that
will be used later.
In particular, it applies to bounded subsets of finite-dimensional,
indeed reflexive, subspaces of $L_{1}$.

\begin{lemma}\label{lem2.1}
Let $H$ be a uniformly integrable subset of  $L_1$ and $\eps > 0$. Then 
there is a $\de > 0$ such that for every  $g \in H$ and every 
$f \in L_1$ with $d(f,0) < \de$ the following inequality holds:
$$
\|f + g\| \ge \|f\| +\| g\| - \eps.
$$
\end{lemma} 

\begin{proof}
Using the uniform integrability of $H$
one can find a $\de > 0$ such that 
$$
2\int_A|g|\,d\mu + 2\de < \eps
$$
for every  $g \in H$ and every measurable
subset $A \subset \Omega$ with $\mu(A) < \de$. 
Now fix
$f \in L_1$ with $d(f,0) < \de$ and denote $A = \{t\dopu  |f(t)| \ge \de\}$.
Then   for every  $g \in H$ we have
\bea
\|f + g\| 
&=& 
\int_A|f + g|\,d\mu   + \int_{\Omega \setminus A}|f + g|\,d\mu  \\
&\ge&  
\int_A|f|\,d\mu   + \int_{\Omega \setminus A} |g|\,d\mu  - 
\int_A|g|\,d\mu   - \int_{\Omega \setminus A} |f|\,d\mu   \\
&\ge& 
\|f\| +\| g\| - 2\int_A|g|\,d\mu  - 2 \de \\
&\ge& 
\|f\| +\| g\| - \eps,
\eea
as requested.
\end{proof}

We now quote Lemma~5.26 from \cite{BL1} that contains the key to the
Bourgain-Rosenthal construction.

\begin{lemma}\label{lem2.2}
Let $0<\eps<1 $. 
Then there is a function 
$f \in L_1[0,1]$ with the following properties:
\begsta
\item 
\label{i*} $f \ge 0$, $\|f\|=1$ and $\|f - \eins\| \ge 2 - \eps$.
\item 
Let $\{f_j\}_{j=1}^{\infty} \subset L_1[0,1]$ be a sequence 
of independent 
random variables with the same distribution as $f$. Then for
every $g \in \Lin\{f_j\}_{j=1}^{\infty}$ with $\|g\| \le 1$
there is a constant function $c$ with $d(g,c) \le \eps$.
\item 
$\|n^{-1} \sum_{j=1}^{n} f_{j} - \eins\| \to 0$ as $n \to \infty$.
\endsta
\end{lemma}

\begin{rem} \label{rem2.2}
Property~(a) of the previous lemma implies in particular
that $d(f,0) \le \sqrt{\eps}$.
Indeed, for 
$A = \{t\dopu  f(t) \ge \sqrt{\eps}\}$ we have ($\lambda$ denotes the
Lebesgue measure)
$$
\int_A|f(t) - 1|\,dt  \le \int_A(f(t)+1)\,dt  - \sqrt{\eps}\lambda(A)
$$
since $(a+1)-|a-1| = 2\min \{a,1\}\ge \sqrt{\eps}$  for
$a\ge\sqrt{\eps}$;
so
$$
 2 - \eps \le  \|f - \eins\|  \le 2 - \sqrt{\eps} \lambda(A)
$$
and hence $\lambda(A)\le \sqrt{\eps}$. 
\end{rem}

In the lemma and in the construction below $\Om$ will be the product 
of countably many copies
of the probability space $[0,1]$. A subspace of $L_1 = L_1\Om$ is
said to depend on finitely many coordinates if all the elements of the 
subspace are functions depending only on a common finite set of
coordinates.

The next lemma is a refinement of \cite[Lemma~5.27]{BL1}; the
difference is that the latter lemma claims (a) only for $u=u_{k}$.

\begin{lemma}\label{lem2.3}
Let $G$ be a finite-dimensional subspace of  $L_1$ that
depends on finitely many coordinates.  
Let $\{u_k\}_{k=1}^{m} \subset S(G)$ and $\eps > 0$. 
Then there is a finite-dimensional subspace $F \subset  L_1$,
also depending on finitely many coordinates and containing $G$,
and there are an integer $n$ and functions 
$\{v_{k,j}\}_{k \le m, j \le n} \subset S(F)$
such that:
\begsta
\item \label{i} 
$\|u + v_{k,j} \| \ge 2 - \eps$ for every $u \in S(G)$ and all $k$
and $j$.
\item \label{ii}
$\|u_k - n^{-1} \sum_{j=1}^{n} v_{k,j} \| \le \eps$ for all $k$.
\item \label{iii} 
For every $\varphi \in B(F)$ 
there is a $\psi \in B(G)$ with $d(\varphi ,\psi ) \le \eps$.
\endsta
\end{lemma}

\begin{proof}
We shall first recall the construction of the proof in \cite[Lemma~5.27]{BL1}
and then point out the necessary changes for our proof.

Let $\eps_1 > 0$ be small enough, $f$ be the function given by
Lemma~\ref{lem2.2} for this $\eps_1$, and $n$ be such that 
$\|n^{-1} \sum_{j=1}^{n} f_{j} - \eins\| \le \eps_1$ for $f_j$ as in 
Lemma~\ref{lem2.2}. (It will become apparent at the end of the proof
how small $\eps_{1}$ should be chosen.)
Let $G$ depend on the first $N$ coordinates of $\Omega$. For 
every $k \le m$ choose $\{f_{k,j}\}_{j \le n}$ which depend on
the
$(N+k)$-th  coordinate of $\Omega$ and are equidistributed with 
$\{f_{j}\}_{j \le n}$. Put $ v_{k,j} = f_{k,j}u_k $, and let 
$F$ be the span of $G$ and of  $\{v_{k,j}\}_{k \le m, j \le n}$.
Note that $\|f_{k,j}u_{k}\|=1$ since the two functions are
stochastically independent. 

The properties (b) and (c) are the same as in 
Lemma~5.27 from \cite{BL1} (cf.\ \cite[p.~118]{BL1} for the norm-1
part of (c)), so we need not repeat their proofs
here. We only have to deal with
property~(a). According to Remark~\ref{rem2.2},  
$d(f_{k,j},0) \le \sqrt{\eps_1}$. Denote 
$$
A  =A_{k,j} = \{t\dopu  f_{k,j}(t) \ge \sqrt{\eps_1}\}, \quad
v_{k,j}^1= f_{k,j}u_k \chi_A, \quad
v_{k,j}^2= f_{k,j}u_k \chi_{\Omega \setminus A}.
$$
Then  
$$
v_{k,j}= v_{k,j}^1 + v_{k,j}^2,
$$
where the first summand has a small support, viz.\  
$$
\mu (\supp v_{k,j}^1) \le \mu (A) \le \sqrt{\eps_1},
$$
 and the second summand has a small
norm, viz.\ 
$$
\|v_{k,j}^2 \| \le \int_{\Omega\setminus A} \sqrt{\eps_{1}} |u_{k}|\,d\mu
\le \sqrt{\eps_1}.
$$ 
So for every $u \in S(G)$ 
$$
\|u + v_{k,j} \| \ge  \|u + v_{k,j}^1\| - \sqrt{\eps_1},
$$
and 
$d( v_{k,j}^1,0) \le \sqrt{\eps_1}$.
To finish the proof
it is enough to apply Lemma~\ref{lem2.1}.
\end{proof}

We now turn to the actual construction of the example.
Fix a decreasing sequence  of numbers $\eps_N > 0$ with 
$\sum_{j=N+1}^\infty \eps_j < \eps_N$ for all $N \in \N$  and
select inductively finite-dimensional subspaces of  $L_1$
$$
\lin\{\eins\}= E_1 \subset E_2 \subset E_3 \subset \dots,
$$  
each of them depending on finitely many coordinates,
$\eps_N$-nets $\{u_k^N\}_{k=1}^{m(N)}$ of $S(E_N)$ and 
collections of elements
$\{v_{k,j}^N\}_{k \le m(N), \, j \le n(N)} \subset S(E_{N+1})$
in such a way that the conclusion of Lemma~\ref{lem2.3} holds with
$\eps=\eps_{N}$, $G=E_{N}$, $F=E_{N+1}$, $\{u_{k}\}_{k=1}^m = 
\{u_{k}^N\}_{k=1}^{m(N)} $, $\{v_{k,j} \}_{k\le m, j\le n}= 
\{v_{k,j}^N\}_{k\le m(N), j\le n(N)} $.
Denote $E = \overline{ \bigcup_{N=1}^\infty E_N }$.

Only part~(a) in the next theorem is new; the other parts are due to
Bourgain and Rosenthal, but we sketch the proofs for completeness.

\begin{thm}\label{th2.1}
The space $E$ constructed above has the following
properties: 
\begsta
\item 
$E$ has the Daugavet property.
\item 
For every $f \in B(E)$ and for every $N \in \N$ 
there exists a $g \in B(E_N)$ such that $d(f,g) < \eps_N$.
\item 
$E$ has the Schur property.
\endsta
\end{thm}
                         
\begin{proof}
(a) According to Lemma~\ref{lem1.2} we need to show that 
for every $u\in S(E)$ and $\eps>0$ the set $\conv(l^+(u,\eps))$ is dense
in $B(E)$. By a perturbation argument it is enough to check
this condition only for  $u$ from the dense
subset $ S(\bigcup_{N=1}^\infty E_N )$ of $S(E)$. 

Fix an $N \in \N$, $u \in S( E_N )$ and $\eps>0$. There is an 
$M > N$ such that $\eps_M < \eps$. By construction 
(property~(a) of Lemma~\ref{lem2.3}) all the elements 
$v_{k,j}^L$ with $L > M$ belong to $l^+(u,\eps)$.
Taking into account property~(b) of Lemma~\ref{lem2.3}
and the fact that  
$\{u_k^L\}_{k=1}^{m(L)}$ forms an $\eps_L$-net of $S(E_L)$,
one can easily establish the density of $\conv(l^+(u,\eps))$ in $B(E)$.

(b) Fix an $f \in B(\bigcup_{J=1}^\infty E_J)$ and an  $N \in \N$. 
Then $f \in B( E_{N+L})$ for some $L$. Applying 
property~(c) of Lemma~\ref{lem2.3} to $f$ we find an 
$f_1 \in S( E_{N+L-1})$ with $d(f,f_1) \le \eps_{N+L}$. 
Applying again
property~(c) of Lemma~\ref{lem2.3} to $f_1$ we find an 
$f_2 \in S( E_{N+L-2})$ with $d(f_1,f_2) \le \eps_{N+L-1}$. 
Continuing in this fashion we obtain
in the $L$-th step some
$g = f_{L}\in S( E_{N}) $ for which 
$$
d(f,g) \le \eps_{N+L}+ \eps_{N+L-1}+ \dots + \eps_{N+1} < \eps_N.
$$

(c) It follows from (b) that the unit ball of $E$ is a precompact
in the metric $d$ of convergence in measure. This 
implies the Schur property.  (In fact, a stronger conclusion can be
drawn, namely, $E$ has the strong Schur property; see
\cite{BourRos}.)
\end{proof}

We observe that the space $E$ cannot be a rich subspace of $L_{1}$ in the
terminology of \cite{KSW}. Indeed, 
 the unit ball of $E$ is precompact in the metric of
convergence in measure, hence its $\tau$-closure $C_E$  is  compact
for that topology~$\tau$. But if $E$ were rich, 
 by \cite[Prop.~2.2]{KKW} $C_E$ would contain $\frac12 B(L_1)$, and
 $B(L_1)$ would be $\tau$-compact as well, which is clearly false.
This remark reveals that $E$ is an essentially new specimen
among the spaces with the \DP.

Since a Banach space with the Schur property cannot contain
infinite-dimensional reflexive subspaces, the following corollary
holds. It answers Question~(4) from \cite{KSW}.

\begin{cor} \label{cor2.6}
There exists a Banach space with the \DP\ that fails to contain
infinite-dimensional reflexive subspaces; in particular it fails to
contain a copy of~$\ell_{2}$.
\end{cor}

It follows as well that the space $E$ from Theorem~\ref{th2.1}
does not contain a copy of $L_{1}$; the
first space with the \DP\ having this feature was constructed in
\cite{KadSSW} after an example given by Talagrand.

One can likewise express Corollary~\ref{cor2.6}  in terms of narrow
operators, a notion studied in detail in  \cite{KSW}. Namely,
the identity operator on $E$ is an operator that does not fix
a copy of $\ell_{2}$, yet it is not narrow. This remark provides a
negative answer to Question~(5) from \cite{KSW}.

Schmidt \cite{Schm} proved that every Dunford-Pettis operator $T$ on          
$L_{1}$ (i.e., $T$ maps weakly convergent sequences to norm
convergent sequences) satisfies the Daugavet equation (\ref{eq1.1}),
in fact, such an operator is easily seen to be narrow on $L_{1}$.
However, on the Schur space $E$ above $-\Id$ is Dunford-Pettis, but
it clearly fails (\ref{eq1.1}). Therefore we have:

\begin{cor} \label{cor2.7}
There is a Banach space with the \DP\ and a Dunford-Pettis operator
on that space which fails {\rm(\ref{eq1.1})}. Hence, in general
Dunford-Pettis operators are not narrow.
\end{cor}


\section{The uniform Daugavet property and ultraproducts}

In Section~6 of the paper \cite{BKSW} conditions for an 
ultraproduct of Banach spaces to possess the Daugavet property
were discussed. For the reader's convenience we recall the relevant
notions and results. 

For a subset $A \subset X$ we denote by $\conv_n(A)$ the set of all convex 
combinations of all $n$-point collections of elements from $A$.
Clearly,
$\conv(A)= \bigcup_{n \in \N} \conv_n(A)$. Denote  
$$
\Daug_n(X,\eps)= \sup_{x,y\in S(X)} \dist(\conv_n(l^+(x,\eps)), y),
$$
where  $\dist(A,B)$ denotes the distance between 
two subsets, i.e., $\dist(A,B)= \inf\{\|a-b\| \dopu a\in A$, $b\in B\}$.

It is easy to see that for every $\eps>0$ the sequence
$(\Daug_n(X,\eps))$ decreases.  

\begin{definition}
A Banach space $X$ is said to have the \textit{uniform Daugavet property} 
if for every $\eps>0$ the sequence $(\Daug_n(X,\eps))$ tends
to 0 when $n$ tends to infinity.
\end{definition}

It follows from Lemma~\ref{lem1.2} that
the uniform Daugavet property implies the Daugavet property.

We now recall a special case of \cite[Th.~6.2]{BKSW}.

\begin{thm} \label{ultraproduct}
Let $\cal U$ be a free ultrafilter on $\N$, $\{X_{i}\}_{i\in \N}$
be a collection of Banach spaces and $X$ be the corresponding
ultraproduct of the $X_{i}$. Then the following assertions are
equivalent:
\begaeq
\item \label{A}
$X$ has the Daugavet property.
\item\label{B}
For every $\eps>0$, $\lim_{\calU,n}\Daug_n(X_{i},\eps)=0$. In
other words, for every fixed $\eps>0$ and  every
$\delta>0$ there is  an $n \in \N$ such
that the set of all $i$ for which $\Daug_n(X_{i},\eps) \le \delta$ 
belongs to the  ultrafilter $\cal U$.                     
\endaeq
In particular, a Banach space $Y$ has the uniform Daugavet property
if and only if  the ultrapower $Y^{\cal U}$ has the Daugavet property.
\end{thm}

\begin{proof}
To deduce (i) from (ii) one just has to notice that if 
$$
\{i\in \N\dopu \Daug_n(X_i,\eps) \le \delta \} \in \calU,
$$ 
then $\Daug_n(X,\eps) \le \delta$.
So  (ii) implies that $\Daug_n(X,\eps)$ tends
to $0$ when $n\to \infty$ for every $\eps>0$, 
which proves the Daugavet property
for $X$.

To deduce (ii) from (i) let us argue ad absurdum. Suppose there
are  $\eps>0$ and  $\delta>0$ such that for every $n \in \N$ the set
$A_n=\{i \in \N\dopu  \Daug_n(X_{i},\eps) > \delta\}$ belongs 
to the  ultrafilter $\cal U$.  The sequence $(A_{n})$ is descending,
but its intersection need not be empty. However,
deleting if necessary a finite number of elements in each of the $A_n$
we can find sets $B_n \in \cal U $, $B_n \subset A_n$, 
$B_1 \supset B_2 \supset \dots $ such that
$\bigcap_n B_n = \emptyset$.

If $i\in B_{n}\setminus B_{n+1}$ we can find $x_{i}, y_{i} \in
S(X_{i})$ such that 
$$
\dist \bigl( \conv_{n}(l^+(x_{i} ,\eps)), y_{i} \bigr) > \delta.
$$
Hence, whenever $i\in B_{n}= \bigcup_{m=n}^\infty B_{m} \setminus
B_{m+1}$ (since $\bigcap_{m}B_{m}= \emptyset$), then
$$
\dist \bigl( \conv_{n}(l^+(x_{i} ,\eps)), y_{i} \bigr) > \delta.
$$
This implies for $x=(x_{i}) \in S(X)$, $y=(y_{i}) \in S(X)$ that
$$
\dist \bigl( \conv_{n}(l^+(x ,\eps)), y \bigr) \ge \delta
$$
for every $n\in \N$ and thus
$$
\dist \bigl( \conv(l^+(x ,\eps)), y     \bigr) \ge \delta
$$
contradicting the \DP\ of $X$.
\end{proof}

It was asked in \cite{BKSW} whether the \DP\ and its uniform version
are actually equivalent. That the Bourgain-Rosenthal spaces may shed
some light on this question is suggested by the following facts. 
The Bourgain-Rosenthal spaces were constructed in order to provide an
example of a Banach space which fails the Radon-Nikod\'ym property
(i.e., some uniformly bounded martingale diverges), yet every
uniformly bounded dyadic martingale converges. Using more
horticultural language, one can express this by saying that the unit
ball contains some $\eta$-bush, but no $\eta$-trees (see Chapter~5 in
\cite{BL1} for these concepts) or indeed no $\eta$-bushes with a
fixed number of branches at each branching node. 

This is reminiscent
of the \DP\ and its uniform variant.  The \DP\ means that  every $y$
of norm~$1$ is almost a convex combination of vectors from
$l^+(x,\eps)$ for any given $x$ of norm~$1$; this enables one to find
an $\eta$-bush for any $\eta<2$. By contrast, the uniform \DP\ is
related to finding such bushes with a fixed number of branches at
each level.

Precisely, we shall now prove the following theorem.

\begin{thm} \label{th3.2} \mbox{ } 
\begsta
\item \label{A1}
For every $n\in \N$ there is a Banach space $X_n$
with the Daugavet property such that 
\begin{equation} \label{En}
\Daug_n \Bigl( X_n,\frac{1}{4} \Bigr) \ge \frac{1}{2}.
\end{equation}
\item\label{B1}
There is a Banach space $X$ which has the Daugavet property,
but does not have the uniform Daugavet property.
\item\label{C1}
An ultrapower of a space with  the Daugavet property
does not necessarily possess the Daugavet property.
\endsta
\end{thm}

\begin{proof}
(a) We take as $X_n$ the  space $E$ 
from Theorem~\ref{th2.1} with parameters $\{\eps_j\}_{j \in \N}$,
where $\eps_1$ is selected in such a way that for every constant
function
$g \in [-2,2]$ and every $f \in L_1$ with $d(f,0) < n \eps_1$ the
inequality 
\begin{equation} \label{bbb}
\|g + f\| \ge \|g\| + \| f\| - \frac{1}{4}
\end{equation}
holds (see Lemma~\ref{lem2.1}). To prove \kref{En} let us check that
\begin{equation} \label{dist}
\dist \bigl(\conv_n(l^+(x,{\textstyle \frac{1}{4}})), y \bigr) \ge \frac{1}{2}
\end{equation}
for the constant functions $x = -\eins $ and $y = \eins $.

Consider an arbitrary element $z \in
\conv_n\bigl(l^+(x,\frac{1}{4})\bigr)$,
$z = \sum_{k=1}^n \lambda_k z_k$, where $\lambda_k \ge 0$,
$\sum_{k=1}^n \lambda_k = 1$, $\|z_k\| \le 1$ and
\begin{equation} \label{zk}
 \|z_k - \eins\| \ge 7/4.
\end{equation}
According to (b) of Theorem~\ref{th2.1} (with $N=1$ 
and $f = z_k$), for every $k \le n$ there is a constant function
$\alpha_k \in [-1,1]$ such that $d(z_k - \alpha_k, 0) =
d(z_{k},\alpha_{k}) < \eps_1$.
Then, using \kref{zk} and \kref{bbb} we conclude
\bea
1 \ge \|z_k\| 
&=& 
\|\alpha_k + (z_k - \alpha_k)\| \\
&\ge& 
|\alpha_k| + \|z_k - \alpha_k\| - \frac{1}{4} \\
&\ge& 
|\alpha_k| +  \|z_k - \eins\| - |1 - \alpha_k|  - \frac{1}{4} \\
&\ge&
|\alpha_k|- |1 - \alpha_k| + \frac32,
\eea
therefore
$$
|\alpha_k|- |1 - \alpha_k| \le - \frac12 .
$$
This implies that $\alpha_k \le 1/4$, consequently 
\begin{equation} \label{viertel}
\sum_{k=1}^n \lambda_k \alpha_k \le 1/4.
\end{equation}
Since 
$d\bigl(\sum_{k=1}^n \lambda_k (z_k - \alpha_k), 0 \bigr) < n \eps_1$,
using \kref{viertel} and \kref{bbb} we deduce that
\bea
\|y-z\| &=& 
\biggl\| \eins - \sum_{k=1}^n \lambda_k z_k \biggr\| \\
&=&
\biggl\| \biggl( \eins - \sum_{k=1}^n \lambda_k \alpha_k \biggr) 
+ \sum_{k=1}^n \lambda_k (\alpha_k - z_k) \biggr\| \\
&\ge&
\biggl\| \eins - \sum_{k=1}^n \lambda_k \alpha_k \biggr\| 
- \frac{1}{4} \ge \frac{1}{2}, 
\eea
which proves \kref{dist}.

(b) It is enough to take the $\ell_{1}$-direct sum $X=X_{1} \oplus_{1}
X_{2} \oplus_{1} \dotsc$; $X$ has the \DP\ by \cite{Woj92}.

(c) This follows from (b) and Theorem~\ref{ultraproduct}.
\end{proof}


\end{document}